\documentclass[11pt]{amsart}

\usepackage{amsthm, amsfonts, amssymb, amscd}
\usepackage[pagebackref,colorlinks]{hyperref}
\usepackage{tikz-cd}
\usepackage{geometry}
\usepackage{marginnote}

\theoremstyle{definition}
\newtheorem{ntn}{Notation}[section]

\theoremstyle{plain}
\newtheorem{lem}[ntn]{Lemma}
\newtheorem{prp}[ntn]{Proposition}
\newtheorem{thm}[ntn]{Theorem}
\newtheorem{cor}[ntn]{Corollary}

\theoremstyle{definition}
\newtheorem{question}[ntn]{Question}
\newtheorem{rem}[ntn]{Remark}
\newtheorem{exa}[ntn]{Example}

\numberwithin{equation}{section}

\newcommand{\N}{\mathbb{N}}
\newcommand{\z}{\mathbb{Z}}
\newcommand{\q}{\mathbb{Q}}
\newcommand{\R}{\mathbb{R}}

\newcommand{\F}{\mathbb{F}}

\newcommand{\EE}{\mathcal{E}}
\newcommand{\Aa}{\mathcal{A}}
\newcommand{\BB}{\mathcal{B}}
\newcommand{\GG}{\mathcal{G}}
\newcommand{\QQ}{\mathcal{Q}}
\newcommand{\WW}{\mathcal{W}}

\newcommand{\PP}{\mathcal{P}}
\newcommand{\RB}{\mathcal{RB}}
\newcommand{\OO}{\mathcal{O}}
\newcommand{\KK}{\mathcal{K}}

\newcommand{\ii}{\mathfrak{i}}

\newcommand{\mmm}{\mathfrak{m}}

\renewcommand{\aa}{{A^\times}}

\newcommand{\tors}{{{\rm Tor}_1^{\z}}}
\newcommand{\exts }{{{\rm Ext}_{\z}^1}}

\newcommand{\lan}{\langle}
\newcommand{\ran}{\rangle}
\newcommand{\se}{\subseteq}
\newcommand{\arr}{\rightarrow}
\newcommand{\larr}{\longrightarrow}
\newcommand{\harr}{\hookrightarrow}
\newcommand{\two}{\twoheadrightarrow}

\newcommand{\PT}{\mathit{{\rm PT}}}

\newcommand{\SL}{\mathit{{\rm SL}}}
\newcommand{\PSL}{\mathit{{\rm PSL}}}

\newcommand{\Hom}{{\rm Hom}}

\renewcommand{\char}{{\rm char}}

\newcommand{\coker}{{\rm coker}}

\newcommand{\ind}{{\rm ind}}
\newcommand{\sign}{{\rm sign}}

\newcommand{\rank}{{\rm rank}}

\newcommand {\mtxx}[4]
{\left(\!\!
	\begin{array}{cc}
		\!\!#1 & \!\!#2   \\
		\!\!#3 & \!\!#4
	\end{array}\!\!
	\right)
}

\newtheoremstyle{athm}
{}
{}
{\itshape}
{}
{\scshape}
{}
{.5em}
{\thmnote{#3}}
\theoremstyle{athm}

\begin{document}
	\title[Low dimensional homology groups of \texorpdfstring{$\SL_2$}{Lg} and 
	\texorpdfstring{$\PSL_2$}{Lg}]
	{On the connections between the low dimensional homology groups of \texorpdfstring{$\SL_2$}{Lg} and 
		\texorpdfstring{$\PSL_2$}{Lg}}
	\author{Behrooz Mirzaii}
	\author{Elvis Torres P\'erez}
	\address{\sf Instituto de Ci\^encias Matem\'aticas e de Computa\c{c}\~ao (ICMC), 
		Universidade de S\~ao Paulo, S\~ao Carlos, Brasil}
	\email{bmirzaii@icmc.usp.br}
	\email{elvistorres@usp.br}
	\begin{abstract}
		In this article we study the low dimensional homology groups of the special linear group
		$\SL_2(A)$ and the projective special linear group $\PSL_2(A)$, $A$ a domain, through the 
		natural surjective map $\SL_2(A) \arr \PSL_2(A)$. In particular, we study the connection 
		of the first, the second and the third homology groups of these groups over euclidean 
		domains $\z[\frac{1}{m}]$, $m$ a square free integer, and local domains. 
	\end{abstract}
	\maketitle
	
	Let $A$ be a commutative ring with unit and $\mu_2(A)$ be the group of $2$-roots of unity in $A$.
	By studying the Lyndon/Hochschild-Serre spectral sequence associated 
	to the central extension
	\[
	1 \arr \mu_2(A) \arr \SL_2(A) \arr \PSL_2(A)\arr 1,
	\]
	and using the homology of cyclic groups \cite[pages 58-59]{brown1994} and the universal coefficient 
	theorem \cite[Exercise 3, \S 1, Chap III]{brown1994} it follows that, for any non-negative integer 
	$n$, the natural homomorphism of homology groups
	\[
	H_n(\SL_2(A),\z) \arr H_n(\PSL_2(A),\z)
	\]
	has 2-power torsion kernel and cokernel. In this article 
	we study  this map for $n\leq 3$.
	
	As our first main theorem we prove that if $A$ is a domain, then we have the 9-term exact sequence
	\[
	0 \arr \mu_2(A)^\sim \arr \widetilde{H}_3(\SL_2(A), \z) \arr H_3(\PSL_2(A),\z)
	\arr \mu_2(A)\otimes_\z H_1(\SL_2(A),\z) 
	\]
	\[
	\arr H_2(\SL_2(A),\z) \arr H_2(\PSL_2(A),\z) \arr \mu_2(A) \arr H_1(\SL_2(A),\z)
	\]
	\[
	\arr H_1(\PSL_2(A),\z) \arr 0,
	\]
	where $\widetilde{H}_3(\SL_2(A), \z)$ is a quotient of $H_3(\SL_2(A),\z)$ over a certain
	2-torsion subgroup and 
	\[
	\mu_2(A)^\sim\simeq \begin{cases}
		0 & \text{ if $\char(A)=2$}\\
		\z/4 & \text{ if $\char(A)\neq 2$}
	\end{cases}
	\]
	(see Theorem \ref{SL-PSL-0} for the precise statement).
	Moreover, if $\char(A)\neq 2$, we show that the image of the map $\mu_2(A)^\sim \arr \widetilde{H}_3(\SL_2(A), \z)$
	coincides with the image of the map $H_3(\lan w\ran,\z) \arr \widetilde{H}_3(\SL_2(A), \z)$, where
	$w={\mtxx{0}{1}{-1}{0}}\in \SL_2(A)$. Observe that the five term exact sequence from the right 
	of the above exact sequence is the well-known five term exact sequence coming from the above central extension
	\cite[Corollary 6.4, Chap. VII]{brown1994}.
	
	As first particular case, we study the above 9-term exact sequence over the euclidean domain $A_m:=\z[\frac{1}{m}]$,
	where $m$ is a square free integer (Proposition \ref{Am}). In particular we show that if $m$ is even, then we have the 
	exact sequence
	\[
	\begin{array}{c}
		0 \arr \z/4 \arr {\displaystyle\frac{H_3(\SL_2(A_m),\z)}{\mu_2(A_m)\otimes_\z H_2(\SL_2(A_m),\z)}}
		\arr H_3(\PSL_2(A_m),\z) \arr 0
	\end{array}
	\]
	and if $m$ is odd, then the sequence
	\[
	\begin{array}{c}
		0 \arr \z/4 \arr {\displaystyle\frac{H_3(\SL_2(A_m),\z)}{\mu_2(A_m)\otimes_\z H_2(\SL_2(A_m),\z)}}
		\arr H_3(\PSL_2(A_m),\z) \arr \z/2 \arr 0
	\end{array}
	\]
	is exact.
	
	As second particular case we study the above 9-term exact sequence over local domains (Proposition 
	\ref{SL-PSL--1}). In particular, if the residue field of a local domain $A$ is different than $\F_2$,
	then we obtain the exact sequence
	\[
	0 \arr \z/4 \arr \frac{H_3(\SL_2(A),\z)}{\mu_2(A)\otimes_\z H_2(\SL_2(A),\z)} \arr H_3(\PSL_2(A),\z) 
	\arr 0.
	\]
	We will see that this is not necessarily true if the residue field of the local domain is $\F_2$. For example we show 
	that for $\z_{(2)}=\{a/b\in \q: a, b \in \z, 2 \nmid b\}$, we have the exact sequence
	\[
	\begin{array}{c}
		0 \arr \z/4 \arr {\displaystyle\frac{H_3(\SL_2(\z_{(2)}),\z)}{\mu_2(\z_{(2)})\otimes_\z H_2(\SL_2(\z_{(2)}),\z)}}
		\arr H_3(\PSL_2(\z_{(2)}),\z) \arr \z/2 \arr 0.
	\end{array}
	\]
	
	In all the above exact sequences the map 
	\[
	\rho_\ast:\mu_2(A)\otimes_\z H_2(\SL_2(A),\z) \arr H_3(\SL_2(A),\z)
	\]
	is induced by the product map 
	\[
	\rho:\mu_2(A) \times \SL_2(A)\arr \SL_2(A), \ \ \ \ (a,X)\mapsto aX.
	\]
	We will show that for $A_2=\z[\frac{1}{2}]$, $\rho_\ast$ is trivial and for $A_3=\z[\frac{1}{3}]$ the map $\rho_\ast$ 
	is non-trivial (see Example \ref{m=2,3}).
	
	
	If $A$ is a local domain such that its residue field has more than two elements, then the natural product map
	\[
	\rho_\ast:\mu_2(A)\otimes_\z H_i(\SL_2(A),\z) \arr H_{i+1}(\SL_2(A),\z)
	\]
	is trivial for $i=0,1$. The following question ask if this is true for $i=2$.
	
	\begin{question}\label{Q1}
		Let $A$ be a local domain such that $|A/\mmm_A|\neq 2$. Is the product map
		\[
		\rho_\ast:\mu_2(A)\otimes_\z H_2(\SL_2(A),\z) \arr H_3(\SL_2(A),\z)
		\]
		trivial?
	\end{question}
	We do not know the answer of the above question even over a general infinite field. But we show that $\rho_\ast$
	is trivial for finite fields, quadratically closed fields and real closed fields (Proposition \ref{some}).
	
	Here we outline the organization of the present paper. In Section 1, we  recall some results  from the literature 
	over the low dimensional homology groups of central extensions. In Section 2 we introduce the scissors congruence 
	group, the Bloch group and formulate the Bloch-Wigner exact sequence over fields. In Section 3 we prove our 9-term 
	exact sequence (Theorem \ref{SL-PSL-0}). In Section 4, we study the above $9$-term exact sequence over the euclidean 
	domain $\z[\frac{1}{m}]$. Finally in Section 5 we investigate the above $9$-term exact sequence over local rings.\\
	~\\
	{\bf Notations.}
	In this paper all rings are commutative, except possibly group rings, and have the unit element $1$.
	For a commutative ring $A$, $\mu(A)$ denotes the group of roots of unity in $A$, i.e.
	\[
	\mu(A):=\{a\in A: \text{there is $n\in \N$ such that $a^n=1$} \},
	\]
	and $\mu_2(A):=\{a\in A: a^2=1\}$. If $\BB \arr \Aa$ is a homomorphism of abelian groups, by $\Aa/\BB$ 
	we mean $\coker(\BB \arr \Aa)$.
	
	\section{The low dimensional homology of central extensions}
	
	Let $\Aa$ be a central subgroup of a group $\GG$. The product map $\rho: \Aa\times \GG \arr \GG$, given by
	$(a,g)\mapsto ag$, induces the natural map $\rho_\ast: H_3(\Aa \times \GG,\z) \arr H_3(\GG,\z)$. By the 
	K\"unneth formula we have the exact sequence
	\[
	\begin{array}{c}
		0 \arr \bigoplus_{i=0}^3 H_i(\Aa,\z)\otimes_\z H_{3-i}(\GG,\z) \arr H_3(\Aa \times \GG,\z) \arr 
		\tors(\Aa, H_1(\GG,\z)) \arr 0.
	\end{array}
	\]
	If 
	\[
	\widehat{H}_3(\Aa \times \GG,\z):=\frac{H_3(\Aa \times \GG,\z)}{\Aa\otimes_\z H_2(\GG,\z)\
		\oplus H_2(\Aa,\z)\otimes_\z H_1(\GG,\z)},
	\]
	then from the above exact sequence we obtain the exact sequence
	\[
	0 \arr H_3(\Aa,\z)\oplus H_3(\GG,\z) \arr \widehat{H}_3(\Aa \times \GG,\z) \arr 
	\tors(\Aa, H_1(\GG,\z)) \arr 0.
	\]
	This exact sequence splits naturally, thus we have a natural splitting
	\[
	\tors(\Aa, H_1(\GG,\z)) \harr \widehat{H}_3(\Aa \times \GG,\z).
	\]
	It is straightforward to check that the image of the composite
	\[
	H_2(\Aa,\z)\otimes_\z H_{1}(\GG,\z)\arr H_3(\Aa \times \GG,\z) \overset{\rho_\ast}{\larr} H_3(\GG,\z)
	\]
	is in the image of the composite
	\[
	\Aa\otimes_\z H_2(\GG,\z)\arr H_3(\Aa \times \GG,\z) \overset{\rho_\ast}{\larr} H_3(\GG,\z)
	\]
	(see \cite[Proposition 4.4(i), Chap. V]{stam1973}). Thus we get the natural composite map
	\[
	\tors(\Aa, H_1(\GG,\z)) \arr \widehat{H}_3(\Aa \times \GG,\z) \arr H_3(\GG,\z)/\rho_\ast(\Aa\otimes_\z H_2(\GG,\z)).
	\]
	We set
	\[
	\widetilde{H}_3(\GG,\z):=H_3(\GG, \z)/\rho_\ast\bigg(\Aa\otimes_\z H_2(\GG,\z)\oplus\tors(\Aa, H_1(\GG,\z)\bigg).
	\]
	
	In this article, for calculations purposes, we will use the (left) bar resolution $B_\bullet(\GG)\arr \z$ 
	of $\z$ over $\z[\GG]$ to study the integral homology groups $H_n(\GG,\z)$ (see \cite[Chap.I, \S 5]{brown1994}):
	\[
	H_n(\GG,\z):=H_n(B_\bullet(\GG)\otimes_{\z[\GG]} \z)=H_n(B_\bullet(\GG)_\GG).
	\]
	If any pair of the elements $g_1, \dots,g_n\in G$ commute, then we define
	\[
	{\pmb c}(g_1,\dots, g_n):= \overline{\sum_{\sigma\in S_n}\sign(\sigma)[g_{\sigma(1)}|\dots|g_{\sigma(n)}]}\in H_n(\GG,\z),
	\]
	where $S_n$ is the symmetric group of degree $n$. If $\Aa:=\lan g_1,\dots, g_n \ran\se \GG$, then 
	the Pontryagin product (see \cite[Chap. V, \S5]{brown1994}) induces a natural map
	\[
	\begin{array}{c}
		\bigwedge_\z^n\Aa \arr H_n(\Aa,\z)\arr H_n(\GG,\z).
	\end{array}
	\]
	It is straightforward to show that the above composite, maps $g_1\wedge \cdots \wedge g_n$ to  ${\pmb c}(g_1,\dots, g_n)$.
	
	A function $\psi: \Aa \arr \BB$ of (additive) abelian groups is called a quadratic map if
	\par (1) for any $a \in \Aa$, $\psi(a)=\psi(-a)$,
	\par (2) the function $\Aa \times \Aa \arr \BB$ with $(a,b) \mapsto \psi(a+b)-\psi(a)-\psi(b)$ is bilinear.\\
	\\
	For any abelian group $\Aa$, there is a universal quadratic map 
	\[
	\gamma: \Aa \arr \Gamma(\Aa)
	\]
	such that for any quadratic map $\psi: \Aa \arr \BB$, there is a unique group homomorphism 
	$\Psi: \Gamma(\Aa) \arr \BB$ such that $\Psi\circ \gamma=\psi$. It is easy to see that $\Gamma$ 
	is a functor from the category of abelian groups to itself. This functor is called Whitehead's quadratic functor.
	In fact if $\EE$ is the free abelian group generated by the symbols $w(a)$, $a \in \Aa$, then we can define 
	$\Gamma(\Aa):=\EE/\KK$, where $\KK$ is the subgroup generated by the relations
	\par (a) $w(a)-w(-a)$=0,
	\par (b) $w(a+b+c)-w(a+b)-w(a+c)-w(b+c)+w(a)+w(b)+w(c)=0$.\\
	The universal quadratic map $\gamma: \Aa \arr \Gamma(\Aa)$ is given by $a\mapsto \gamma(a):=w(a)+\KK$.
	For more on Whitehead's quadratic functor we refer the reader to \cite[\S1.2]{baues1996}.
	
	Now let $\Aa$ be a central subgroup of $\GG$. The inclusion $i:\Aa\arr \GG$ induces
	the natural map $i_\ast: \Aa\arr H_1(\GG,\z)$, $a\mapsto {\pmb c}(a)$. 
	By the universal property of $\Gamma$, the quadratic map
	\[
	\psi:\Aa \arr \Aa\otimes_\z H_1(\GG,\z),  \ \ \  a \mapsto a\otimes {\pmb c}(a),
	\]
	can be extended to a group homomorphism
	\begin{equation}\label{tau}
		\tau: \Gamma(\Aa) \arr \Aa\otimes_\z H_1(\GG,\z).
	\end{equation}
	Since for $\rho_\ast: a \otimes {\pmb c}(g) \mapsto {\pmb c}(a,g)$, the composite
	$A \overset{\psi}{\larr} \Aa\otimes_\z H_1(\GG,\z) \overset{\rho_\ast}{\larr} H_2(\GG,\z)$
	is trivial, the composite
	\[
	\Gamma(\Aa) \overset{\tau}{\arr} \Aa\otimes_\z H_1(\GG,\z) \overset{\rho_\ast}{\larr} H_2(\GG,\z)
	\]
	is trivial too. This induces the natural map 
	\[
	\coker(\tau )\arr H_2(\GG,\z).
	\]
	Let $\GG$ be a group. Let $K(\GG, n)$ be an Eilenberg-MacLane space of type $(\GG,n)$. This is 
	a CW-complex such that
	\[
	\pi_k(K(\GG, n), \ast)\simeq
	\begin{cases}
		1 & \text{if $k\neq n$} \\ 
		\GG & \text{if $k = n$}.
	\end{cases}
	\]
	This space is unique up to homotopy. Note that for $n\geq 2$, $\GG$ should be abelian. 
	It is well-known \cite[Theorem 21.1]{EM1954} that
	\[
	\Gamma(\Aa)\simeq H_4(K(\Aa, 2),\z).
	\]
	
	Let $1\arr \Aa \arr \GG \arr  \QQ\arr 1$ be a central extension. The standard classifying space 
	theory gives a (homotopy  theoretic) fibration \cite[Chap. 6]{dk2001} of Eilenberg-MacLane spaces 
	\[
	K(\Aa, 1) \arr K(\GG, 1)\arr K(\QQ, 1)
	\]
	(see \cite[5.1.28]{r1996}). From this we obtain the fibration 
	\[
	K(\GG, 1)\arr K(\QQ, 1) \arr K(\Aa, 2)
	\]
	(see \cite[Lemma 3.4.2]{mp2012}). From the Serre spectral sequence of this fibration we obtain 
	the sequence
	\[
	H_4(K(\QQ,1),\z) \arr H_4(K(\Aa, 2),\z) \overset{\tau}{\larr} H_2(K(\Aa,2),\z)\otimes_\z H_1(K(\GG,1), \z).
	\]
	This gives us the sequence $H_4(\QQ,\z) \arr \Gamma(\Aa) \overset{\tau}{\larr} \Aa\otimes_\z H_1(\GG, \z)$,
	which from it we obtain the natural map 
	\[
	H_4(\QQ,\z) \arr \ker(\tau).
	\]
	The following theorem is proved by Eckmann and Hilton in \cite{eh1971}, which generalises Ganea's 6-term
	exact sequence \cite[Theorem 2.2, Chap. V]{stam1973}.
	
	\begin{thm}[Eckmann-Hilton]
		\label{EH}
		For any central extension $1\arr \Aa \arr \GG \arr  \QQ \arr 1$, we have the natural $10$-term exact sequence
		\[
		H_4(\QQ, \z) \arr \ker(\tau) \arr \widetilde{H}_3(\GG,\z) \arr H_3(\QQ,\z)
		\arr \coker(\tau )\arr H_2(\GG,\z) \arr  H_2(\QQ,\z) 
		\]
		\[
		\arr \Aa \arr H_1(\GG,\z)\arr H_1(\QQ,\z) \arr 0.
		\]
	\end{thm}
	\begin{proof}
		See \cite[Theorem~1.1]{eh1971}. 
	\end{proof}
	
	Let $\Aa$ be a finite cyclic group. If $2\mid |\Aa|$, let $\Aa^\sim$ denote the unique non-trivial 
	extension of $\Aa$ by $\z/2$. If $2\nmid|\Aa|$, we define $\Aa^\sim:=\Aa$. Thus if $n=|\Aa|$, then
	\[
	\Aa^\sim \simeq \begin{cases}
		\z/2n & \text{if $2\mid n$}\\
		\z/n & \text{if $2\nmid n$.}
	\end{cases}
	\]
	
	Let $\Aa$ be an
	ind-cyclic group, i.e. $\Aa$ is direct limit of its finite cyclic subgroups. Then by passing 
	to the limit we can define $\Aa^\sim$. Note that if $\Aa$ has no $2$-torsion, then $\Aa^\sim=\Aa$.
	There is always a canonical injective homomorphism $\Aa \harr \Aa^\sim$ whose composition with the 
	projection $\Aa^\sim \arr \Aa$ coincides with multiplication by $2$. 
	
	\begin{lem}\label{G-I}
		If $\Aa$ is an ind-cyclic group, then $\Gamma(\Aa)\simeq \Aa^\sim$.
	\end{lem} 
	\begin{proof}
		This follows from the facts that the functor $\Gamma$ commutes with direct limit and 
		\[
		\Gamma(\z/n)\simeq \z/{\rm gcd}(n^2,2n)\simeq (\z/n)^\sim,
		\]
		where ${\rm gcd}(n^2,2n)$ is the greatest common divisor of $n^2$ and $2n$ (see \cite[Corollary 1.2.9]{baues1996}).
	\end{proof}
	
	Let $\Aa$ be an abelian group. Let $\sigma_1:\tors(\Aa,\Aa)\arr \tors(\Aa,\Aa)$ be obtained by interchanging 
	the copies of $\Aa$. This map is induced by the involution $\Aa\otimes_\z \Aa \arr \Aa\otimes_\z \Aa$, 
	$a\otimes b\mapsto -b\otimes a$ \cite[\S2]{mmo2022}. Let $\Sigma_2'=\{1, \sigma'\}$ be the symmetric group 
	of order 2 and consider the following action of $\Sigma_2'$ on $\tors(\Aa,\Aa)$:
	\[
	(\sigma', x)\mapsto -\sigma_1(x).
	\]
	
	\begin{prp}\label{H3A}
		Let $\Aa$ be an abelian group and $T_\Aa$ be its torsion subgroup. Then we have the exact sequence
		\[
		\begin{array}{c}
			0 \arr \bigwedge_\z^3 \Aa\arr H_3(\Aa,\z) \arr \tors(T_\Aa,T_\Aa)^{\Sigma_2'}\arr 0.
		\end{array}
		\]
		If $T_\Aa$ is an ind-cyclic group, then $\Sigma_2'$ acts trivially on $\tors(T_\Aa,T_\Aa)$  
		and the exact sequence 
		\[
		\begin{array}{c}
			0 \arr \bigwedge_\z^3 \Aa\arr H_3(\Aa,\z) \arr \tors(T_\Aa,T_\Aa) \arr 0
		\end{array}
		\]
		splits naturally.
	\end{prp}
	\begin{proof}
		By \cite[Lemma 5.5]{suslin1991} or \cite[Section~\S6]{breen1999} we have the exact sequence 
		\[
		\begin{array}{c}
			0 \arr \bigwedge_\z^3 \Aa\arr H_3(\Aa,\z) \arr \tors(\Aa,\Aa)^{\Sigma_2'} \arr 0.
		\end{array}
		\]
		Since $\tors(\Aa,\Aa)\simeq \tors(T_\Aa,T_\Aa)$, we obtain the first exact sequence.
		Now let $T_\Aa$ be an  ind-cyclic group.  Since $\tors(\z/n,\z/n)\simeq \z/n$, the action of $\Sigma_2'$ 
		on $\tors(\z/n,\z/n)$ is trivial. Now by passing to the limit, we see that the action of $\Sigma_2'$ on
		$\tors(T_\Aa,T_\Aa)$ is trivial. For the last part note that since $T_\Aa$ is ind-cyclic, 
		$\bigwedge_\z^3 T_\Aa=0$. Now applying the first part to the inclusion $T_\Aa\harr \Aa$, we obtain 
		the following commutative diagram with exact rows
		\[
		\begin{tikzcd}
			&&H_3(T_\Aa,\z) \ar[r, "\simeq"]\ar[d]& \tors(T_\Aa, T_\Aa)\ar[d, "\simeq"]&\\
			0\ar[r] & \bigwedge_\z^3 \Aa \ar[r] & H_3(\Aa,\z)\ar[r] & \tors(\Aa,\Aa)\ar[r]&0.
		\end{tikzcd}
		\]
		Now from this diagram we obtain a natural splitting map 
		\[
		\tors(T_\Aa, T_\Aa)\simeq H_3(T_\Aa,\z) \arr H_3(\Aa,\z).
		\]
	\end{proof}
	
	\begin{cor}\label{H3A-1}
		For any integer $m\in \z$, let $m:\Aa \arr \Aa$ be given by $a\mapsto ma$. Then the map 
		$m_\ast:H_3(\Aa,\z) \arr H_3(\Aa,\z)$ induces multiplication by $m^3$ on $\bigwedge_\z^3 \Aa$ 
		and multiplication by $m^2$ on $\tors(\Aa,\Aa)^{\Sigma_2'}$.
	\end{cor}
	\begin{proof}
		This follows from the above proposition.
	\end{proof}
	
	\section{The bloch-wigner exact sequence}
	Let $A$ be a commutative ring and set $\WW_A:=\{a\in A: a(1-a)\in \aa\}$.
	The {\it scissors congruence group} $\PP(A)$ of $A$ is defined as the 
	quotient of the free abelian group  generated by symbols $[a]$, $a\in \WW_A$, by the subgroup 
	generated by the elements
	\[
	[a] -[b]+\bigg[\frac{b}{a}\bigg]-\bigg[\frac{1- a^{-1}}{1- b^{-1}}\bigg]+ \bigg[\frac{1-a}{1-b}\bigg],
	\]
	where $a, b, a/b  \in \WW_A$. Let 
	$S_\z^2(\aa):=(\aa \otimes_\z \aa)/\lan a\otimes b+b\otimes a:a,b \in \aa\ran$. The map
	\[
	\lambda: \PP(A) \arr S_\z^2(\aa), \ \ \ \ [a] \mapsto a \otimes (1-a),
	\]
	is well-defined and its kernel is called the {\it Bloch group} of $A$ 
	and is denoted by $\BB(A)$. 
	
	For any local ring  $A$ there is a natural map of graded rings $K_\bullet^M(A) \to K_\bullet(A)$, 
	from Milnor $K$-theory to $K$-theory. It is easy to see that 
	\[
	K_0^M(A)\simeq K_0(A)\simeq \z,  \ \ \ \ K_1^M(A)\simeq K_1(A)\simeq \aa.
	\]
	Moreover,  Matsumoto 
	proved that for a field $F$
	\[
	K_2(F)\simeq K_2^M(F) \simeq (\aa \otimes_\z \aa)/\lan a\otimes (1-a): a\in \WW_F \ran
	\]
	(see \cite[Theorem 4.3.15]{r1996}).
	
	The \textit{indecomposable $K_3^\ind(F)$ of a field $F$} is defined as follows:
	\[
	K_3^\ind(F):
	=K_3(F)/K_3^M(F).
	\]
	The Bloch group and  the indecomposable part of the third $K$-group are deeply connected.
	
	\begin{thm}[Bloch-Wigner exact sequence]\label{bloch-wigner}
		For any field $F$, we have a natural exact sequence 
		\[
		0 \arr \tors(\mu(F),\mu(F))^\sim \arr K_3^\ind (F) \arr \BB(F) \arr 0.
		\]
	\end{thm}
	\begin{proof}
		The case of infinite fields has been proved in \cite[Theorem 5.2]{suslin1991} and the case of finite fields has 
		been settled in \cite[Corollary 7.5]{hutchinson-2013}. 
	\end{proof}
	
	For any filed $F$, there are natural maps $H_2(\SL_2(F),\z) \arr K_2(F)$ and 
	$H_3(\SL_2(F),\z) \arr K_3^\ind(F)$ (see \cite[\S5]{hutchinson-2013}). 
	
	\begin{prp}\label{K2-K3}
		If $F$ is a quadratically closed field, then we have natural isomorphisms 
		\[
		K_2(F)\simeq K_2^M(F)\simeq H_2(\SL_2(F),\z),  \ \ \ \ \ 
		K_3^\ind(F)\simeq H_3(\SL_2(F),\z).
		\]
		Furthermore, $K_2(F)$ is uniquely $2$-divisible.
	\end{prp}
	\begin{proof}
		See \cite[page 286, 2.12]{sah1989}, \cite[Proposition 6.4]{mirzaii-2008} and 
		\cite[Proposition 1.2]{bass-tate1973}.
	\end{proof}
	
	The following classical Bloch-Wigner exact sequence will be used in the proof of Theorem~\ref{SL-PSL-0}.
	
	\begin{cor}[Bloch-Wigner, Dupont-Sah]\label{CBW}
		If $F$ is a quadritically closed field, then we have the exact sequence
		\[
		\begin{tikzcd}
			0\arr \tors(\mu(F),\mu(F))\arr H_3(\SL_2(F),\z)\arr\BB(F)\arr 0.
		\end{tikzcd}
		\]
	\end{cor}
	\begin{proof}
		For a direct proof see \cite[Appedix A]{dupont-sah1982}. See also \cite[Theorem~5.2]{suslin1991}.
	\end{proof}
	
	Note that the map $\mu(F) \arr \SL_2(F)$ given by $a\mapsto {\mtxx{a}{0}{0}{a^{-1}}}$ induces the natural map
	\[
	H_3(\mu(F),\z) \arr H_3(\SL_2(F),\z).
	\]
	By Proposition \ref{H3A}, $H_3(\mu(F),\z)\simeq \tors(\mu(F),\mu(F))$.
	Now the left map in the above classical Bloch-Wigner exact sequence is the map
	\[
	\tors(\mu(F),\mu(F))\simeq H_3(\mu(F),\z) \arr H_3(\SL_2(F),\z).
	\]
	
	\begin{prp}[Coronado \cite{cor2021}]\label{coronado}
		Let $F$ be a real closed field. Then we have the Bloch-Wigner exact sequence
		\[
		\begin{tikzcd}
			0\arr \tors(\mu(F),\mu(F))^\sim \arr H_3(\SL_2(F),\z)\arr\BB(F)\arr 0.
		\end{tikzcd}
		\]
		In particular, $H_3(\SL_2(F), \z)\simeq K_3^\ind(F)$
	\end{prp}
	\begin{proof}
		See \cite[Theorem 6.1]{cor2021} and its proof.
	\end{proof}
	
	\section{Low dimensional homology of \texorpdfstring{$\SL_2$}{Lg} and \texorpdfstring{$\PSL_2$}{Lg}}
	
	We apply Theorem \ref{EH} to the central extension 
	\[
	1 \arr \mu_2(A) \overset{i}{\arr} \SL_2(A) \overset{p}{\arr} \PSL_2(A) \arr 1,
	\]
	where $\mu_2(A)=\{a\in \aa:a^2=1\}$ and the map $i$ is given by $a\mapsto {\mtxx{a}{0}{0}{a}}={\mtxx{a}{0}{0}{a^{-1}}}$.
	Let
	\[
	\widetilde{H}_3(\SL_2(A), \z):=\displaystyle\frac{H_3(\SL_2(A), \z)}{\rho_\ast\Big(\mu_2(A)\otimes_\z H_2(\SL_2(A),\z)\oplus 
		\tors(\mu_2(A), H_1(\SL_2(A),\z))\Big) },
	\]
	where $\rho:\mu_2(A) \times \SL_2(A)\arr \SL_2(A)$ is the product map and is given by $(a, X)\mapsto aX$. 
	We recall that for a field $F$,
	\[
	\mu_2(F)^\sim \simeq \begin{cases}
		0 & \text{if $\char(F)=2$}\\
		\z/4 & \text{if $\char(F)\neq 2$}.
	\end{cases}
	\]
	Now we are ready to state and prove our first main result.
	
	\begin{thm}\label{SL-PSL-0}
		Let $A$ be a ring such that there is a ring homomorphism $A\arr F$, $F$ a field, where 
		$\mu_2(A) \simeq \mu_2(F)$. Then we have the $9$-term exact sequence
		\[
		0 \arr \mu_2(A)^\sim \arr \widetilde{H}_3(\SL_2(A), \z) \arr H_3(\PSL_2(A),\z) 
		\arr \mu_2(A)\otimes_\z H_1(\SL_2(A),\z) 
		\]
		\[
		\arr H_2(\SL_2(A),\z) \arr H_2(\PSL_2(A),\z) \arr \mu_2(A) \arr H_1(\SL_2(A),\z)
		\]
		\[
		\arr H_1(\PSL_2(A),\z) \arr 0.
		\]
		Moreover, if $\char(F)\neq 2$, then the image of the maps $\mu_2(A)^\sim \arr \widetilde{H}_3(\SL_2(A), \z)$ and
		$H_3(\lan w \ran,\z)\arr \widetilde{H}_3(\SL_2(A), \z)$ coincide, where $w:={\mtxx{0}{1}{-1}{0}}$.
	\end{thm}
	\begin{proof}
		First we study the homomophism $\tau:\Gamma(\mu_2(A))\arr \mu_2(A)\otimes_\z H_1(\SL_2(A),\z)$ (see (\ref{tau})). 
		Note that by Lemma \ref{G-I}, $\Gamma(\mu_2(A))\simeq \mu_2(A)^\sim$. The map $\tau$ is induced by the quadratic map
		\[
		\mu_2(A)\arr \mu_2(A)\otimes_\z H_1(\SL_2(A),\z),\ \ \ \ -1\mapsto (-1)\otimes i_\ast(-1),
		\]
		where $i_\ast:\mu_2(A)\simeq H_1(\mu_2(A),\z) \arr H_1(\SL_2(A),\z)$ is given by $-1 \mapsto {\pmb c}(-I_2)$.
		Since $w^2=-I_2$, we have
		\[
		(-1)\otimes i_\ast(-1)=(-1)\otimes {\pmb c}(-I_2)=(-1)\otimes {\pmb c}(w^2)
		=(-1)\otimes 2 {\pmb c}(w)=(-1)^2\otimes {\pmb c}(w)=0.
		\]
		Thus $\tau$ must be trivial. Now by Theorem \ref{EH} we obtain the exact sequence
		\[
		H_4(\PSL_2(A),\z)\arr \Gamma(\mu_2(A))\arr\widetilde{H}_3(\SL_2(A),\z) \arr H_3(\PSL_2(A),\z) 
		\]
		\[
		\arr \mu_2(A)\otimes_\z H_1(\SL_2(A),\z) \arr \cdots.
		\]
		Therefore to complete the proof of the first part of the theorem, we must prove that the map 
		\[
		H_4(\PSL_2(A),\z)\arr \Gamma(\mu_2(A))\simeq \mu_2(A)^\sim
		\]
		is trivial. Since this map is trivial when $\mu_2(A)\simeq\mu_2(F)=1$, we may assume that 
		$\mu_2(A)\simeq\mu_2(F)=\{\pm1\}\simeq \z/2$.
		
		Let $\overline{F}$ be the algebraic closure of $F$. Then $H_1(\SL_2(\overline{F}),\z)=0$
		and $H_2(\SL_2(\overline{F}),\z)$ is uniquely 2-divisible by Proposition \ref{K2-K3}.
		These imply that $\tors(\mu_2(\overline{F}), H_1(\SL_2(\overline{F}),\z))$ and
		$\mu_2(\overline{F})\otimes_\z H_2(\SL_2(\overline{F}),\z)$
		are trivial. Hence by Theorem \ref{EH}, we have the exact sequence
		\begin{equation}\label{inj}
			\Gamma(\mu_2(\overline{F})) \arr H_3(\SL_2(\overline{F}),\z)\arr H_3(\PSL_2(\overline{F}),\z)\arr 0.
		\end{equation}
		Observe that $\Gamma(\mu_2(\overline{F}))\simeq \mu_2(\overline{F})^\sim\simeq \z/4$ (Lemma \ref{G-I}).
		The surjective map $p: \SL_2(\overline{F}) \arr \PSL_2(\overline{F})$ induces the
		commutative diagram with exact rows
		\[
		\begin{tikzcd}
			0\ar[r] &\tors(\mu(\overline{F}),\mu(\overline{F}))\ar[r]\ar[d, "\alpha"]&
			H_3(\SL_2(\overline{F}),\z)\ar[r]\ar[d,"p_\ast"]&\mathcal{\BB}(\overline{F})\ar[r]\ar[d,equal]&0\\
			0 \ar[r] & \tors (\widetilde{\mu}(\overline{F}), \widetilde{\mu}(\overline{F})) \ar[r] & 
			H_3(\PSL_2(\overline{F}),\z) \ar[r] & \BB(\overline{F}) \ar[r] & 0.
		\end{tikzcd}
		\]
		The first row is the classical Bloch-Wigner exact sequence discussed in the previous section 
		(Corollary \ref{CBW}).  The second exact sequence can be proved as the first exact sequence, where 
		$\widetilde{\mu}(\overline{F}):=\mu(\overline{F})/\mu_2(\overline{F})$ (see \cite[Theorem 4.10 and 
		Appendix A]{dupont-sah1982}. Also see \cite[Corollary 3.5]{B-E-2024} below). The left hand side square 
		is obtained from the commutative diagram
		\[
		\begin{tikzcd}
			\tors(\mu(\overline{F}),\mu(\overline{F}))\ar[r,equal]\ar[d, "\alpha"]& H_3(\mu(\overline{F}),\z)\ar[r,hook]
			\ar[d, "p_\ast"]&H_3({\rm T}(\overline{F}),\z)\ar[r]\ar[d, "p_\ast"] &H_3(\SL_2(\overline{F}),\z)\ar[d, "p_\ast"]\\
			\tors(\widetilde{\mu}(\overline{F}),\widetilde{\mu}(\overline{F}))\ar[r,equal]&
			H_3(\widetilde{\mu}(\overline{F}),\z)\ar[r,hook] & H_3(\PT(\overline{F}),\z) \ar[r] & H_3(\PSL_2(\overline{F}),\z),
		\end{tikzcd}
		\]
		where for a ring $R$, 
		\[
		{\rm T}(R):=\bigg\{ {\mtxx{a}{0}{0}{a^{-1}}}|a\in R^\times \bigg\},
		\ \ \ \ \ 
		\PT(R):={\rm T}(R)/\mu_2(R)I_2
		\]
		and the middle square is induced by the map $\mu(R)\arr {\rm T}(R)$, $a\mapsto {\mtxx{a}{0}{0}{a^{-1}}}$.
		The extension  
		\[
		1 \arr \mu_2(\overline{F}) \arr {\rm T}(\overline{F}) \arr \PT(\overline{F}) \arr 1
		\]
		is isomorphic with the extension 
		\[
		1 \arr\mu_2(\overline{F})\arr\overline{F}^\times \overset{()^2}{\larr} \overline{F}^\times\arr 1. 
		\]
		The map $\phi=()^2:\overline{F}^\times \arr \overline{F}^\times$ induces the map
		$\phi_\ast: H_3(\overline{F}^\times,\z) \arr H_3(\overline{F}^\times,\z)$. 
		By Proposition \ref{H3A} we have the split exact sequence
		\[
		\begin{array}{c}
			0 \arr \bigwedge_\z^3\overline{F}^\times \arr H_3(\overline{F}^\times,\z) \arr 
			\tors(\mu(\overline{F}), \mu(\overline{F}))\arr 0.
		\end{array}
		\]
		By Corollary \ref{H3A-1} this leads to multiplication by $2^2$ on 
		$\tors(\mu(\overline{F}),\mu(\overline{F}))$. Since 
		\[
		\tors(\mu(\overline{F}),\mu(\overline{F}))\simeq \mu(\overline{F}),
		\]
		the above fact implies that $\ker(\alpha)\simeq\z/4$. On the other hand by the Snake lemma 
		$\ker(\alpha)\simeq \ker(p_\ast)$. Thus the left map of the exact sequence (\ref{inj}) must 
		be injective and so we have the exact sequence
		\[
		0\arr \Gamma(\mu_2(\overline{F})) \arr H_3(\SL_2(\overline{F}),\z) \arr 
		H_3(\PSL_2(\overline{F}),\z) \arr 0. 
		\]
		Observe that
		\[
		\Gamma(\mu_2(\overline{F}))\simeq \mu_2(\overline{F})^\sim \harr \mu(\overline{F})^\sim \simeq \mu(\overline{F})\simeq 
		\tors(\mu(\overline{F}),\mu(\overline{F})).
		\]
		The composite $A\arr F \harr \overline{F}$ induces the morphism of extensions
		\[
		\begin{tikzcd}
			1\ar[r] &\mu_2(A)\ar[r]\ar[d]&\SL_2(A)\ar[r]\ar[d]& \PSL_2(A)\ar[r]\ar[d]&1\\
			1\ar[r] &\mu_2(\overline{F}) \ar[r] & \SL_2(\overline{F})\ar[r] &\PSL_2(\overline{F})\ar[r]&1.
		\end{tikzcd}
		\]
		Now by Theorem \ref{EH} we obtain the commutative diagram with exact rows
		\[
		\begin{tikzcd}
			H_4(\PSL_2(A),\z)\ar[r]\ar[d]&\Gamma(\mu_2(A)) \ar[r]\ar[d]&\widetilde{H}_3(\SL_2(A),\z) \ar[r]\ar[d] &
			H_3(\PSL_2(A),\z)\ar[d]& \\
			0 \ar[r] & \Gamma(\mu_2(\overline{F})) \ar[r] & H_3(\SL_2(\overline{F}),\z) \ar[r] & 
			H_3(\PSL_2(\overline{F}),\z) \ar[r] & 0,
		\end{tikzcd}
		\]
		where the exactness of the second row has just been shown. Since $\mu_2(A)\simeq \mu_2(F)=\mu_2(\overline{F})$, 
		we have 
		\[
		\Gamma(\mu_2(A))\simeq \Gamma(\mu_2(\overline{F})).
		\]
		Now from the above diagram it follows that
		the map $\Gamma(\mu_2(A)) \arr\widetilde{H}_3(\SL_2(A),\z)$ is injective. This complete the proof of the 
		first claim.
		
		For the second claim, first we show that the composite
		\[
		H_3(\lan w \ran,\z)  \arr \widetilde{H}_3(\SL_2(A), \z)  \arr H_3(\PSL_2(A),\z)
		\]
		is trivial. Let denote this map by $\varphi_\ast$. The group $H_3(\lan w \ran,\z)$ is cyclic of 
		order 4 and is generated by
		\[
		\chi:=\overline{[w|w|w]+[w|\!-\!1|w]+[w|\!-\!w|w]} \in H_3(B_\bullet(\lan w\ran)_{\lan w \ran})
		=H_3(\lan w \ran,\z),
		\]
		(see \cite[Proposition 8.1]{zikert2015} or \cite[Proposition 3.25]{parry-sah1983}),
		where $B_\bullet(\lan w\ran)\arr \z$ is the bar resolution of $\z$ over $\lan w\ran$. We have
		\[
		\varphi_\ast(\chi)=\overline{[w|w|w]+[w|\!-\!w|w]}=\overline{2[w|w|w]}.
		\]
		Since 
		\[
		d_4([w|w|w|w])={2[w|w|w]}
		\]
		in $B_\bullet(\PSL_2(A))_{\PSL_2(A)}$, we have $\varphi_\ast(\chi)=0$.
		This shows that the above composite is trivial. Now by the commutative diagram
		\[
		\begin{tikzcd}
			H_3(\lan w \ran,\z) \ar[r]\ar[d, equal]& \widetilde{H}_3(\SL_2(A),\z) \ar[d]\\
			H_3(\lan w \ran,\z) \ar[r]& H_3(\SL_2(\overline{F}),\z),
		\end{tikzcd}
		\]
		it is sufficient to prove that the map $H_3(\lan w \ran,\z) \arr H_3(\SL_2(\overline{F}),\z)$
		is injective. Let $\ii^2=-1$. Then $w$ is conjugate to ${\mtxx{\ii}{0}{0}{\ii^{-1}}}$. It is known
		by the classical Bloch-Wigner exact sequence (Corollary \ref{CBW}) that the map
		\[
		H_3(\mu(\overline{F}),\z)\arr H_3(\SL_2(\overline{F}),\z)
		\]
		which is induced by $\mu(\overline{F}) \arr \SL_2(\overline{F})$, $a \mapsto {\mtxx{a}{0}{0}{a^{-1}}}$, 
		is injective. Thus the map 
		\[
		H_3(\lan \ii \ran,\z) \arr H_3(\SL_2(\overline{F}),\z)
		\]
		is injective. This proves our claim.
	\end{proof}
	
	\begin{rem}
		The above theorem is also valid if we replace $\SL_2(A)$ with the elementary subgroup ${\rm E}_2(A)$.
		The proof stays the same. Just one should note that over a field $F$, ${\rm E}_2(F)=\SL_2(F)$.
	\end{rem}
	
	\begin{cor}
		Let $A$ be a ring such that there is a ring homomorphism $A\arr F$, $F$ a field of $\char(F)\neq 2$, 
		where $\mu_2(A) \simeq \mu_2(F)$. Then the map $H_3(\lan w \ran,\z)\arr H_3(\SL_2(A), \z)$ is injective.
	\end{cor}
	\begin{proof}
		By Theorem \ref{SL-PSL-0}, the map $H_3(\lan w \ran,\z)\arr \widetilde{H}_3(\SL_2(A), \z)$ is injective.
		But this map factors through $H_3(\SL_2(A), \z)$. Thus we have our claim.
	\end{proof}
	
	
	\section{The low dimensional homology of \texorpdfstring{$\SL_2(\z[\frac{1}{m}])$}{Lg}
		and \texorpdfstring{$\PSL_2(\z[\frac{1}{m}])$}{Lg}}
	
	It is known that  $\SL_2(\z)\simeq \z/4 \ast_{\z/2} \z/6$ and $\PSL_2(\z)\simeq \z/2 \ast \z/3$.
	By the Mayer-Vietoris exact sequence associated to a amalgamated product \cite[Corollary 7.7, Chap. II]{brown1994} 
	one obtains 
	\[
	H_1(\SL_2(\z),\z)\simeq \z/12,\ \ \ \ \ H_2(\SL_2(\z),\z)=0,
	\]
	\[
	H_3(\SL_2(\z),\z)\simeq \z/12, \ \ \ H_3(\PSL_2(\z),\z)\simeq \z/6. 
	\]
	Using these and Theorem \ref{SL-PSL-0} we obtain the exact sequence
	\begin{equation}\label{ZZZ}
		0 \arr \z/4 \arr H_3(\SL_2(\z),\z) \arr H_3(\PSL_2(\z),\z) \arr \z/2 \arr 0,
	\end{equation}
	which will be used in the proof of the following result.
	
	\begin{prp}\label{Am}
		Let $m$ be a square free integer and let $A_m:=\z[\frac{1}{m}]$. 
		\par {\rm (i)} If $m$ is even, then $H_1(\SL_2(A_m),\z)\simeq H_1(\PSL_2(A_m),\z)$ and the sequences 
		\[
		0\arr H_2(\SL_2(A_m),\z)\arr H_2(\PSL_2(A_m),\z) \arr \mu_2(A_m) \arr 0,
		\]
		\[
		\begin{array}{c}
			0 \arr \z/4 \arr {\displaystyle\frac{H_3(\SL_2(A_m),\z)}{\mu_2(A_m)\otimes_\z H_2(\SL_2(A_m),\z)}}
			\arr H_3(\PSL_2(A_m),\z) \arr 0
		\end{array}
		\]
		are exact.
		\par {\rm (ii)} If $m$ is odd, then $H_2(\SL_2(A_m),\z)\simeq H_2(\PSL_2(A_m),\z)$ and the sequences 
		\[
		0\arr \mu_2(A_m) \arr H_1(\SL_2(A_m),\z)\arr H_1(\PSL_2(A_m),\z) \arr  0,
		\]
		\[
		\begin{array}{c}
			0 \arr \z/4 \arr {\displaystyle\frac{H_3(\SL_2(A_m),\z)}{\mu_2(A_m)\otimes_\z H_2(\SL_2(A_m),\z)}}
			\arr H_3(\PSL_2(A_m),\z) \arr \z/2 \arr 0
		\end{array}
		\]
		are exact.
	\end{prp}
	\begin{proof}
		By \cite[Proposition 4.3]{B-E2024} or \cite[Theorem 1.2]{nyberg2024} we have
		\[
		H_1(\SL_2(A_m),\z)\simeq 
		\begin{cases}
			0 & \text{if $2\mid m$, $3\mid m$}\\
			\z/3 & \text{if $2\mid m$, $3\nmid m$}\\
			\z/4 & \text{if $2\nmid m$, $3\mid m$.}\\
			\z/12 & \text{if $2\nmid m$, $3\nmid m$}
		\end{cases}
		\]
		It follows from this that the inclusion $\z \se A_m$ induces the surjective splitting map
		\[
		\z/12\simeq H_1(\SL_2(\z),\z)\two H_1(\SL_2(A_m),\z).
		\]
		Now from the commutative diagram
		\[
		\begin{tikzcd}
			\tors(\mu_2(\z),H_1(\SL_2(\z),\z))  \ar[r, "\rho_\ast=0"]\ar[d, two heads]& H_3(\SL_2(\z),\z)\ar[d]\\
			\tors(\mu_2(A_m),H_1(\SL_2(A_m),\z)) \ar[r, "\rho_\ast"]& 
			\displaystyle\frac{H_3(\SL_2(A_m),\z)}{\mu_2(A_m)\otimes_\z H_2(\SL_2(A_m),\z)}
		\end{tikzcd}
		\]
		it follows that the bottom map is trivial. In the above diagram, first note that since $H_2(\SL_2(\z),\z)=0$, we have
		$H_3(\SL_2(\z),\z)=\displaystyle\frac{H_3(\SL_2(\z),\z)}{\mu_2(\z)\otimes_\z H_2(\SL_2(\z),\z)}$ and thus the above
		diagram is commutative. Second, by the exact sequence (\ref{ZZZ}), the upper map 
		in the above diagram should be trivial.
		
		(i) Let $m$ be even. Then from the above calculation of $H_1(\SL_2(A_m),\z)$, we see that
		\[
		\mu_2(A_m)\otimes H_1(\SL_2(A_m),\z)=0
		\]
		and the natural map $\mu_2(A_m) \arr H_1(\SL_2(A_m),\z)$
		is trivial. Thus the claim follows from Theorem \ref{SL-PSL-0}.
		
		(ii) Now assume that $m$ is odd. Then $\mu_2(A_m)\otimes H_1(\SL_2(A_m),\z)\simeq \z/2$.
		Moreover, from the commutative  diagram
		\[
		\begin{tikzcd}
			H_3(\PSL_2(\z),\z) \ar[r, two heads]\ar[d]& \mu_2(\z)\otimes_\z H_1(\SL_2(\z),\z) \ar[d, two heads]\\
			H_3(\PSL_2(A_m),\z) \ar[r]& \mu_2(A_m)\otimes_\z H_1(\SL_2(A_m),\z)
		\end{tikzcd}
		\]
		it follows that the bottom map is surjective (note that the upper map is surjective in view of 
		Theorem \ref{SL-PSL-0}, because $H_2(\SL_2(\z),\z)=0$). Since 
		\[
		H_1(\PSL_2(A_m),\z)\simeq (\z/2)^\alpha \oplus (\z/3)^\beta
		\]
		for some $\alpha$ and $\beta$ (see \cite[Corollary 4.4]{ww1998}), from the group structure of 
		$H_1(\SL_2(A_m),\z)$ one sees that the natural map 
		\[
		\mu_2(A_m)\arr H_1(\SL_2(A_m),\z)
		\]
		must be injective. Hence the claim follows from Theorem \ref{SL-PSL-0}.
	\end{proof}
	
	\begin{rem}
		A central extension $1\arr \Aa \arr \GG \arr  \QQ \arr 1$ is called a stem-extension if the natural map 
		\[
		\Aa = H_1(\Aa,\z)\arr H_1(\GG,\z)=\GG^{\rm ab}
		\]
		is trivial \cite[\S4]{ehs1972} and it is  called a weak stem-extension if the natural map 
		\[
		\Aa \otimes_\z \Aa \arr \Aa \otimes_\z H_1(\GG,\z)
		\]
		is trivial \cite[page 103]{ehs1972}. The proof of Theorem \ref{SL-PSL-0} shows that if $A$ is a ring 
		such that $\mu_2(A)=\{\pm1\}$, then  
		\[
		1 \arr \mu_2(A) \arr \SL_2(A) \arr \PSL_2(A) \arr 1
		\]
		always is a weak stem-extension. But by Proposition \ref{Am} this extension is not a stem-extension
		for $A=\z[\frac{1}{m}]$, where $m$ is a square free odd integer. This also is true for $A=\z$ by
		the paragraph above Proposition \ref{Am} and $A=\z_{(2)}$ by Example \ref{Z2}.
	\end{rem}
	
	We learned the following result from Kevin Hutchinson. It seems that this lemma is widely known, 
	but we could not find any suitable reference for it.
	
	\begin{lem}\label{coho-ho}
		Let $\GG$ be a group such that $H_{n}(\GG,\z)$, $H_{n+1}(\GG,\z)$
		and $H^{n+1}(\GG,\z)$ are finitely generated. Then 
		\[
		\rank\ H_{n+1}(\GG,\z)=\rank\ H^{n+1}(\GG,\z),\ \ \ and\ \ \ H_n(\GG,\z)_{\rm tor}\simeq H^{n+1}(\GG,\z)_{\rm tor}.
		\]
	\end{lem}
	\begin{proof}
		This claim easily follows from the Universal Coefficients Theorem 
		\[
		0\arr \exts(H_n(\GG,\z),\z) \arr H^{n+1}(\GG,\z) \arr \Hom(H_{n+1}(\GG,\z),\z)\arr 0
		\]
		and the facts that for any abelian group $\Aa$, $\exts(\z,\Aa)=0$ and $\exts(\z/n,\Aa)\simeq \Aa/n\Aa$.
	\end{proof}

	\begin{exa}\label{m=2,3}
		Let $A_m:=\z[\frac{1}{m}]$, where $m=2$ or $3$. 
		\par (i)  The product map $\mu_2(A_2)\otimes_\z H_2(\SL_2(A_2),\z) \arr H_3(\SL_2(A_2),\z)$ is trivial:
		By \cite[Corolllary 3.7, Remark 3.8]{B-E-2024} we have the isomorphism 
		\[
		H_3(\PSL_2(A_2),\z)\simeq \RB(A_2),
		\]
		where $\RB(A_2)$ is the refined Bloch group of $A_2$ (see \cite[\S4]{C-H2022} for the definition of $\RB(A_2)$). 
		By a result of Coronado and Hutchinson, we have the refined Bloch-Wigner exact sequence
		\[
		0\arr \z/4 \arr H_3(\SL_2(A_2),\z) \arr \RB(A_2) \arr 0
		\]
		(see \cite[Proposition 8.31]{C-H2022}). From these two results and Proposition \ref{Am}, it follows that the 
		desired product map is trivial.
		\par (ii) Now we show that the product map $\mu_2(A_3)\otimes_\z H_2(\SL_2(A_3),\z) \arr H_3(\SL_2(A_3),\z)$ 
		is not trivial. We know that 
		\[
		H_3(\SL_2(A_3), \z)\simeq\z/12 \oplus \z/4.
		\]
		This follows from the cohomology calculation of $\SL_2(A_3)$ in \cite[page 8]{adem-naffah1998}
		and Lemma \ref{coho-ho}. Furthermore, by  \cite[Proposition 4.2]{ww1998}, $H^k(\PSL_2(A_3), \z)$ 
		is generated by elements of orders $2$ and $3$ for $k\geq 3$. Again by Lemma \ref{coho-ho} 
		it follows that $H_3(\PSL_2(A_3), \z)$ is generated by elements of orders $2$ and $3$. Now 
		let the above product map be trivial. Then by Proposition~\ref{Am} the sequence
		\[
		\begin{array}{c}
			0 \arr \z/4 \arr H_3(\SL_2(A_3),\z) \arr H_3(\PSL_2(A_3),\z) \arr \z/2 \arr 0
		\end{array}
		\]
		must be exact.
		This and the isomorphism $H_3(\SL_2(A_3),\z)\simeq \z/12 \oplus \z/4$ imply that $H_3(\PSL_2(A_3),\z)$
		has an element of order $4$, which is a contradiction. Observe that $H_2(\SL_2(A_3),\z)\simeq \z$ and thus
		$\mu_2(A_3)\otimes_\z H_2(\SL_2(A_3),\z)\simeq \z/2$.
	\end{exa}
	
	\section{The low dimensional homology groups of \texorpdfstring{$\SL_2$}{Lg}
		and \texorpdfstring{$\PSL_2$}{Lg} of local rings}
	
	Let $A$ be a local ring with maximal ideal $\mmm_A$. It is known that 
	\[
	H_1(\SL_2(A),\z)\simeq
	\begin{cases}
		A/\mmm_A^2 &  \text{if $|A/\mmm_A|=2$}  \\
		A/\mmm_A &  \text{if $|A/\mmm_A|=3$.}  \\
		0 &  \text{if $|A/\mmm_A|\geq  4$}  \\
	\end{cases}
	\]
	(see \cite[Proposition 4.1]{B-E2024}). Now let $|A/\mmm_A|\neq 2$.
	Since $\mu_2(A)\otimes_\z H_1(\SL_2(A),\z)=0$, from Theorem \ref{EH} we obtain the exact sequence
	\[
	0 \arr H_2(\SL_2(A),\z) \arr  H_2(\PSL_2(A),\z) \arr \mu_2(A) \arr 0.
	\]
	Since
	\[
	\mu_2(A)\otimes_\z H_1(\SL_2(A),\z)=0,  \ \ \ \ \tors(\mu_2(A), H_1(\SL_2(A),\z))=0,
	\]
	by Theorem \ref{SL-PSL-0} we obtain the following result.
	
	\begin{prp}\label{SL-PSL--1}
		Let $A$ be a local ring such that there is a ring homomorphism $A\arr F$, $F$ a field, where 
		$\mu_2(A) \simeq \mu_2(F)$. If $|A/\mmm_A|\neq 2$, then the sequence 
		\[
		0 \arr \mu_2(A)^\sim \arr\frac{H_3(\SL_2(A),\z)}{\mu_2(A)\otimes_\z H_2(\SL_2(A),\z)} 
		\arr H_3(\PSL_2(A),\z) \arr 0
		\]
		is exact.
	\end{prp}
	%
	
	\begin{rem}
		The result of the above proposition is valid for any domain $A$ such that $H_1(\SL_2(A),\z)$ is 
		torsion of odd exponent. For example let $\OO_{-3}$  be the ring of algebraic integers of the 
		imaginary quadratic field $\q(\sqrt{-3})$. The homology groups of $\SL_2(\OO_{-3})$ and 
		$\PSL_2(\OO_{-3})$ have been calculated in \cite[Theorem, page 376]{alperin1980} and 
		\cite[Theorem 5.7]{S-V1983}. We list these homology groups for $1\leq n\leq 3$:
		\[
		H_n(\SL_2(\OO_{-3}),\z)\simeq 
		\begin{cases}
			\z/3 & \text{if $n=1$} \\  
			\z/4 & \text{if $n=2,$} \\
			\z/24 \times \z/6 & \text{if $n=3$} 
		\end{cases}
		\]
		\[
		H_n(\PSL_2(\OO_{-3}),\z)\simeq 
		\begin{cases}
			\z/3 & \text{if $n=1$} \\  
			\z/4 \times \z/2& \text{if $n=2.$} \\
			\z/6 \times \z/3 & \text{if $n=3$} 
		\end{cases}
		\]
		From Theorem \ref{SL-PSL-0} we have the exact sequence
		\[
		\begin{array}{c}
			0 \arr \z/4 \arr {\displaystyle\frac{H_3(\SL_2(\OO_{-3}),\z)}{\mu_2(\OO_{-3})\otimes_\z 
					H_2(\SL_2(\OO_{-3}),\z)}} \arr H_3(\PSL_2(\OO_{-3}),\z) \arr 0.
		\end{array}
		\]
		From the sizes of the third homology groups of $\SL_2(\OO_{-3})$ and $\PSL_2(\OO_{-3})$ 
		it follows that the product map
		\[
		\mu_2(\OO_{-3})\otimes_\z H_2(\SL_2(\OO_{-3}),\z) \arr H_3(\SL_2(\OO_{-3}),\z)
		\]
		is non-trivial.
	\end{rem}
	
	\begin{exa}\label{Z2}
		Let $A$ be a local domain such that $A/\mmm_A\simeq \F_2$. Then $H_1(\SL_2(A),\z)\simeq A/\mmm_A^2$.
		By a careful analysis of the composition
		\[
		\mu_2(A) \arr H_1(\SL_2(A),\z)\overset{\simeq}{\larr} A/\mmm_A^2,
		\]
		one sees that $-1\in \mu_2(A)$ maps to $-\overline{6}\in A/\mmm_A^2$. Since $2\in \mmm_A$, 
		we have $-\overline{6}=\overline{2}$ in $A/\mmm_A^2$. Thus if $2\in \mmm_A^2$, then 
		$H_1(\SL_2(A),\z)\simeq H_1(\PSL_2(A),\z)$. Otherwise we have the exact sequence
		\[
		1\arr \mu_2(A) \arr H_1(\SL_2(A),\z) \arr H_1(\PSL_2(A),\z)\arr 0.
		\]
		Observe that
		\[
		\mu_2(A)\otimes_\z H_1(\SL_2(A),\z)\simeq\mu_2(A)\otimes_\z A/\mmm_A^2\simeq A/(2A+\mmm_A^2).
		\]
		If $A=\z_{(2)}=\{a/b\in \q: a,b \in \z, 2\nmid b\}$ (or $A=\z_2$, the ring of $2$-adic 
		integers), then $2\notin \mmm_A^2$ (by the Nakayama lemma). Thus the sequence
		\[
		1\arr \mu_2(\z_{(2)}) \arr H_1(\SL_2(\z_{(2)}),\z) \arr H_1(\PSL_2(\z_{(2)}),\z)\arr 0
		\]
		is exact. Note that $H_1(\SL_2(\z_{(2)}),\z)\simeq \z/4$ \cite[Example 4.2]{B-E2024}.
		With an argument similar to the proof of Proposition \ref{Am}(ii), one can show that  
		\[
		\tors(\mu_2(\z_{(2)}),H_1(\SL_2(\z_{(2)}),\z)) \arr
		\displaystyle\frac{H_3(\SL_2(\z_{(2)}),\z)}{\mu_2(\z_{(2)})\otimes_\z H_2(\SL_2(\z_{(2)}),\z)}
		\]
		is trivial and 
		\[
		H_3(\PSL_2(\z_{(2)}),\z) \arr \mu_2(\z_{(2)})\otimes_\z H_1(\SL_2(\z_{(2)}),\z)\simeq \z/2
		\]
		is surjective. Thus 
		\[
		H_2(\SL_2(\z_{(2)}),\z)\simeq H_2(\PSL_2(\z_{(2)}),\z)
		\]
		and the sequence 
		\[
		\begin{array}{c}
			0 \arr \z/4 \arr {\displaystyle\frac{H_3(\SL_2(\z_{(2)}),\z)}{\mu_2(\z_{(2)})\otimes_\z 
					H_2(\SL_2(\z_{(2)}),\z)}} \arr H_3(\PSL_2(\z_{(2)}),\z) \arr \z/2 \arr 0
		\end{array}
		\]
		is exact. We have similar results for $\z_2$. Hence Proposition \ref{SL-PSL--1} is not necessarily true
		when $|A/\mmm_A|=2$.
	\end{exa}
	
	\begin{prp}\label{some}
		Let $F$ be either a finite field, a quadratically closed field, or a real closed field. Then the sequence
		\[
		0\arr \mu_2(F)^\sim \arr H_3(\SL_2(F),\z)\arr  H_3(\PSL_2(F),\z) \arr 0
		\]
		is exact. In particular, the product map 
		\[
		\mu_2(F)\otimes_\z H_2(\SL_2(F),\z) \arr H_3(\SL_2(F),\z)
		\]
		is trivial.
	\end{prp}
	\begin{proof}
		(i) Let $F=\F_q$ be finite. For fields of characteristic $2$ the above exact sequence is trivial. 
		Thus we may assume that $q$ is odd. We have seen that 
		\[
		H_1(\SL_2(\F_q),\z)\simeq
		\begin{cases}
			0 & \text{if $q\neq 3$}\\
			\z/3 & \text{if $q= 3$.}
		\end{cases}
		\]
		Moreover, it is well-known that
		\[
		H_2(\SL_2(\F_q),\z)\simeq
		\begin{cases}
			0 & \text{if $q\neq 9$}\\
			\z/3 & \text{if $q= 9$}
		\end{cases}
		\]
		(see \cite[page 280]{sah1989}). These facts together with Theorem \ref{SL-PSL-0} give us the desired 
		exact sequence.
		
		(ii) Let $F$ be  quadratically closed. By Proposition \ref{K2-K3}, $H_2(\SL_2(F), \z)$ is 
		uniquely 2-divisible. Thus by Proposition \ref{SL-PSL--1} we obtain the desired exact sequence.
		
		(iii) Let $F$ be a real closed field. By Proposition \ref{coronado} we have the 
		Bloch-Wigner exact sequence
		\[
		0\arr \tors(\mu(F),\mu(F))^\sim \arr H_3(\SL_2(F),\z)\arr\BB(F)\arr 0.
		\]
		From \cite[Corollary 3.5]{B-E-2024} and \cite[Corollary 4.7]{cor2021} we obtain the isomorphism 
		\[
		H_3(\PSL_2(F),\z)\simeq\BB(F).
		\]
		Our claim follows from these two facts. Thus we presented a proof of this in above proposition.
	\end{proof}
	
	\begin{rem}
		It is mentioned in \cite[page 230]{dps1988} that the sequence
		\[
		0 \arr \z/4 \arr H_3(\SL_2(\R),\z) \arr H_3(\PSL_2(\R),\z) \arr 0
		\]
		is exact and for a proof the authors referred to \cite[App. A]{dupont-sah1982}. We couldn't really 
		verify this. So we gave a proof of this in above proposition.
	\end{rem}
	

\end{document}